\documentclass[12pt]{article}

\usepackage{amsthm,amsfonts,amssymb,latexsym}

\def\F{\mathcal{F}}
\def\A{\overline{A}}
\def\P{\mathcal{P}}
\def\L{\mathcal{L}}
\def\R{\mathcal{R}}
\def\I{\mathcal{I}}
\def\S{\mathcal{S}}
\def\Aut{\mathop{\mbox{\rm Aut}}\nolimits}
\def\IS{\mathcal{I}\mathcal{S}}

\newtheorem{lem}{Lemma}[section]
\newtheorem{prp}{Proposition}[section]
\newtheorem{thm}{Theorem}[section]
\newtheorem{cor}{Corollary}[section]
\newtheorem{definition}{Definition}[section]
\newtheorem{example}{Example}[section]
\newtheorem{rmk}{Remark}[section]

\title{
Decomposition theorem for invertible substitutions on three-letter
alphabet
\thanks{Research supported by NSFC and by the Special Funds for Major
State Basic Research Projects of China}}

\author
{Bo TAN\thanks{Department of Mathematics, Wuhan University,
430072 Hubei, Wuhan, P. R. China,
e-mail: \ tanbo@colmath.whu.edu.cn, zhxwen@whu.edu.cn, ypzhang@whu.edu.cn
}\and Zhi-Xiong WEN$^\dag$ \and Yiping ZHANG$^\dag$}

\date{}

\begin{document}

\maketitle


\begin{abstract}
We study the structure of invertible substitutions on three-letter
alphabet. We show that there exists a finite set ${\mathbb S}$ of
invertible substitutions such that any invertible substitution can
be written as $I_w\circ \sigma_1\circ\sigma_2\circ\cdots\circ
\sigma_k$, where $I_w$ is the inner automorphism associated
with~$w$, and $\sigma_j\in{\mathbb S}$
 for $1\le j\le k$. As a consequence, $M$ is
the matrix of an invertible substitution if and only if it is a finite
product of non-negative elementary matrices.

\end{abstract}



\section{Introduction}
\setcounter{equation}{0}

The study of substitutions (endomorphisms of of the free monoid of
finite type) plays an important role in finite automata, symbolic
dynamics, and fractal geometry (\cite{book1,AI,Co,Qu,Su}). It has various applications in quasicrystals, computational
complexity, information theory (see for instance
\cite{ABI2,AG,BT1,LGJJ}). In addition, substitution is also a
fundamental object studied in combinatorial group theory
\cite{Lo2,Ly,MKS}.

For substitutions over two-letter alphabet plenty of results have been
obtained \cite{Be2,MS,Se,WWsp}. The notion of invertible substitution
appears in~\cite{PWW3}: these are the substitutions which extend as an
automorphism of the corresponding free group. Since then, they have
been studied by many authors (see for example \cite{Be2,EI,WWss}).

The invertible substitutions over a two-letter alphabet form a monoid
whose structure in known~\cite{WWli}: this monoid can be generated by
a permutation and two so-called Fibonacci substitutions (see Theorem
\ref{thm1.1} below).

This result on the structure of invetible substitutions had many
applications, namely to the study of local isomorphisms of fixed
points of substitutions~\cite{WWli,WWW} and to the study of trace
maps~\cite{PWW3}.


When the alphabet has more than two letters, the situation is much
more complicated. In \cite{WZ}, it was shown, by enumerating
infinitely many so-called indecomposable substitutions, that the monoid of
invertible substitutions over three letters (which will be denoted by
$\IS(A^*)$) is not finitely generated. So, up to now the structure of
$\IS(A^*)$ remained unknown.


In this paper, we elucidate the structure of $\IS(A^*)$.  We show that
if~$\sigma$ is an invertible substitution over a three-letter
alphabet, there exists a word~$w$ such that $I_w\circ \sigma$ or
$I_{w^{-1}}\circ \sigma$ is the composition of finitely many Fibonacci
substitutions and permutations (Theorem \ref{thm2.1}). As a
consequence, the substitution matrix of an invertible substitution is
positively decomposable (Theorem \ref{thm2.2}).

\section{Preliminaries and notations}

Let us first recall some basic definitions and notations in the theory
of substitutions (see \cite{book1,Qu, MKS} for a general theory).

Let $A=\{a, b, c\}$ (resp.\
$\overline{A}=\{a,b,c,a^{-1},b^{-1},c^{-1}\}$), $A^*$ (resp.\
$\Gamma_A$) be the free monoid (resp.\ the free group) generated by
$A$ (the unit element is the empty word $\varepsilon$).  The elements
in $A^*$ will be called ``positive words" or
simply ``words" and those of $\Gamma_A$ ``signed'' or ``mixed''
words. The inverse of a positive word will be said to be ``negative''.

Let $w\in \Gamma_A$. If $w=x_1\cdots x_k$ with $x_i\in \A\
(i=1,2,\cdots,k)$ and $x_ix_{i+1}\not=\varepsilon$ ($i=1,\cdots,k-1$),
we say that $x_1\cdots x_k$ is in the reduced form
and that the length of $w$ is $k$. This length will be denoted
by~$|w|$.


Let $w,w_j\ (j=1,\cdots,k)\in\Gamma_A$. If $w=w_1w_2\cdots w_k$
satisfies $|w|=|w_1|+\cdots+|w_k|$, we say that $w_1w_2\cdots w_k$ is
a reduced expression of $w$. We then say that $w_1$ is a prefix of $w$
and that $w_k$ is a suffix of $w$, and we then write $w_1\lhd w$ and
$w_k\rhd w$ respectively.

As already defined, a substitution over~$A$ is a morphism
$\sigma$ of $A^*$. Such a morphism extends in a natural way to an
endomorphism of $\Gamma_A$; If this extension is an automorphism of
$\Gamma_A$, the substitution~$\sigma$ is said to be invertible. The
set of substitutions (resp.\ invertible substitutions) is denoted by
$\S(A^*)$ (resp.\ by $\IS(A^*)$).

We often identify an endomorphism~$\sigma$ of~$\Gamma_A$ with the
triple $\bigl(\sigma(a),\sigma(b),\sigma(c)\bigr)$ of (maybe
mixed) words.  We define the length of~$\sigma$ to be
$|\sigma|=|\sigma(a)|+|\sigma(b)|+|\sigma(c)|$.

If $U$ is a subset of a monoid, $\langle U\rangle$ stands for the
sub-monoid (not the sub-group, even when dealing inside a group)
generated by $U$. The use of the same notation for different groups
and monoids will not generate any confusion.
\bigskip

We shall use the following basic invertible substitutions and
automorphisms.

\bigskip
\noindent{\bf $\bullet$ Permutations}
\bigskip


$\P$ will denote the symmetric group on $A$.
Notice that $\P=\langle\pi_1,\ \pi_2\rangle$,
where $\pi_1=(b,a,c),\ \pi_2=(c,b,a)$.
Note that the identity $I=(a,b,c)\in \P$.

%

\bigskip
\noindent{\bf $\bullet$ Fibonacci type}:

\bigskip

Set $\phi_l=(ba, b, c)$, $\phi_r=(ab, b, c)$.

Set $\L=\{\pi\circ\phi_l\circ\pi';\ \pi,\pi'\in\P\}$,
$\R=\{\pi\circ\phi_r\circ\pi';\ \pi,\pi'\in\P\}$,
and $\F=\R\bigcup\L$.

The elements in $\F$ are called
substitutions of Fibonacci type or simply Fibonacci substitutions.

\bigskip
\noindent{\bf $\bullet$ Simple substitutions}

\bigskip
It is easy to check the following equalities.
$$
S:=\langle\pi_1, \pi_2, \phi_l,\phi_r\rangle=
\langle
\P, \L, \R
\rangle
\subset\IS(A^*)
$$
\begin{definition}\label{8888}
The elements in $S$ will be called simple substitutions.
\end{definition}

\bigskip
\noindent{\bf $\bullet$ Involutions}

\bigskip
$\iota_1:=(a^{-1}, b, c)$, $\iota_2:=(a, b^{-1}, c)$, $\iota_3:= (a,
b, c^{-1})$, $\I:=\{\iota_i\ ;\ 1\le i\le 3\}$.
\bigskip


It is well known (the cancellation theory of Nielsen, see for example
\cite{MKS}) that
$$\Aut(\Gamma_A)=\langle{\P}, {\mathcal F},{\mathcal I}\rangle
=\langle\pi_1, \pi_2, \phi_l, \phi_r,\iota_1\rangle.$$

\bigskip
For the case $A=\{a,b\}$, the following decomposition is known.

\begin{thm}\label{thm1.1}\emph{\cite{WWli}}\quad
Let $A_2=\{a,b\},\ \alpha=(b,a),\ \beta=(ab,b),\ \gamma=(ba,b)$.  Then
we have
\begin{equation}\label{eqn1.1}
\IS(A_2^*)=\langle\alpha,\beta, \gamma\rangle.
\end{equation}
\end{thm}

\begin{rmk}\label{rmk1.1}
With our terminology, the above theorem says that, for the case of
two-letter alphabet, every invertible
substitution is simple.
In the case $A=\{a,b,c\}$ it is shown in \emph{\cite{WZ}} that
there exist non-simple substitutions
and that there exists no finite set of
substitutions which generates $\IS(A^*)$ . Consequently the above theorem cannot be extended to the
case of three letters. Our Theorem \ref{thm2.1} below will give a
characterization of the structure of $\IS(A^*)$.
\end{rmk}

\bigskip
\noindent{\bf $\bullet$ Indecomposable substitutions}

\bigskip
\begin{definition}\label{dec}

Let $\sigma$ be an invertible substitution.

$\sigma$ is called trivial (resp.\ non-trivial) if $\sigma\in\P$ (resp.\
$\sigma\not\in\P$).

If there exist non-trivial invertible substitutions $\sigma_1$ and $\sigma_2$
such that $\sigma=\sigma_1\circ\sigma_2$, we say that $\sigma$ is
decomposable. Otherwise, we say that $\sigma$ is
indecomposable.

\end{definition}

The elements in ${\P}\cup{\mathcal F}$ are indecomposable
in the sense of this definition. They are ``simple indecomposable
substitutions" according to Definition \ref{8888}.
The others will be called non-simple indecomposable
substitutions. As a matter of fact, there exist infinitely many indecomposable
substitutions, hence $\IS(A^*)$ is not finitely generated~\cite{WZ}.

\begin{example}\label{ex1.1}
Here are some examples of non-simple indecomposable substitutions.

\bigskip
$
\begin{array}{lr}
  (ab, acb,ac^n)\ (n\ge 2) & \\ (ab^{n+1}c,ab^nc,acc)\ (n\ge 1)& \\
  (acbc^n,ac^{n+1},ab^kc^n)\ (n,k\ge 1) &\\
\end{array}
$
\end{example}

\bigskip
Notice that any invertible substitution can be written as the product
of finitely many indecomposable substitutions.

\bigskip
\noindent{\bf $\bullet$ Inner automorphisms}

\bigskip
Let $z\in \Gamma_A$. $I_z\in\Aut(\Gamma_A)$ is defined as follows:
$I_z(w)=zwz^{-1}\ (w\in \Gamma_A).$ That is $I_z=(zaz^{-1},zbz^{-1},
zcz^{-1})$. We have $I_z\circ(w_1,w_2,w_3)=(zw_1z^{-1},zw_2z^{-1},
zw_3z^{-1})$.  Notice also that $I_{\varepsilon}=I$.

\bigskip
We shall also use the following operation.

\bigskip
\noindent{\bf $\bullet$ Cyclic operators:}

\bigskip
Let $w=x_1x_2\cdots x_k\in A^*$. Then for $1\le i\le k$,
$$
\begin{array}{rl}
I_{(x_1\cdots x_i)^{-1}}(w):& = x_{i+1}\cdots x_kx_1\cdots x_i \\
I_{x_{i+1}\cdots x_n}(w):& = x_{i+1}\cdots x_kx_1\cdots x_i
\end{array}
$$
We observe that the above operator is in fact a $k$-cyclic
operator.

The following operations on substitutions
will be frequently used in this paper:
$$
\begin{array}{rl}
  I_{z^{-1}}\circ(zw_1,zw_2,zw_3) & =(w_1z,w_2z,w_3z), \\
  I_{z}\circ(w_1z,w_2z,w_3z) & =(zw_1,zw_2,zw_3).
\end{array}
$$
More generally, let $\sigma\in\S(A^*)$, $z\in A^*$ (or $z^{-1}\in A^*)$
be such that
$\sigma'=I_z\circ \sigma$ is still a substitution,
we shall say that $I_z$
is a cyclic operator for $\sigma$.

\bigskip
\noindent{\bf $\bullet$ Substitution matrix:}

\bigskip
Let $w\in \Gamma_A$.  $|w|_a$ (resp. $|w|_b$, $|w|_c$) will denote the
number (algebraic sum) of appearances of $a$ (resp. $b$, $c$) in $w$
(e.g. $|a^{-1}ba|_a=0$).

Let $\sigma\in \S(A^*)$. The substitution matrix of $\sigma$ is defined
by:
$$M_\sigma=\pmatrix{ |\sigma(a)|_a & |\sigma(b)|_a & |\sigma(c)|_a \cr
|\sigma(a)|_b & |\sigma(b)|_b & |\sigma(c)|_b \cr |\sigma(a)|_c &
|\sigma(b)|_c & |\sigma(c)|_c }$$

Let $\sigma ,\tau\in \S(A^*)$. We have
$M_{\sigma\circ\tau}=M_{\sigma}M_{\tau}$.  If $\sigma\in \IS(A^*)$,
then $\det(M_\sigma)=\pm 1$.

The above definitions and equalities
can be extended to the case where $\sigma$ and $\tau$ are endomorphisms
on $\Gamma_A$.

Substitution matrix is a basic tool to study substitutions \cite{Qu},
but it cannot characterize the non-commutative combinatorial
properties of $A^*$ and $\S(A^*)$.

Notice that the usual non-negative elementary matrices are
exactly the substitution matrices of permutations or Fibonacci
substitutions.

As have been pointed out in \cite{WZ}, comparing the case of two
letters with that of three letters, we have

\begin{rmk}\label{rmk1.2}

1) If $M=(m_{ij}),\ (m_{ij}\in \{0,1,2,\cdots\})$ is a $2\times
2$-matrix with determinant $\det(M)=\pm 1$, then $M$ is a finite
product of non-negative elementary matrices. In particular, any such
matrix is the substitution matrix of some invertible substitution.

\bigskip
2) The above statement is no longer true for $3\times 3$-matrices.
\end{rmk}
A simple counter example is \cite{WZ}
$$M=\pmatrix{ 3 & 0 & 1 \cr 0 & 2 & 1 \cr 1 & 1 & 1}. $$ We can verify
that $\det(M)=1$, that $M$ can not be decomposed as a product of
non-negative elementary matrices, and that it is not a substitution
matrix of any invertible substitution .

We mention that there are infinitely many such matrices.

\bigskip
It was an open problem whether the substitution matrix of an
invertible substitution is a finite product of non-negative elementary
matrices \cite{book1}.  Our Theorem \ref{thm2.2} will give an affirmative answer.

\section{Main theorems}
\setcounter{equation}{0}

The
following decomposition theorem characterizes the structure of
$\IS(A^*)$: any invertible substitution is a simple one up to a
cyclic operator.

\begin{thm}\label{thm2.1}
  Let $\sigma\in\IS(A^*)$. There exists $w\in A^*$ or $w^{-1}\in A^*$
such that

(1) \ $I_w\circ \sigma$ is a simple substitution. In other words,
after a cyclic operation ($I_w\circ \sigma$), any invertible
substitution becomes simple.

(2) Furthermore, we can take $w$ (or $w^{-1}$) to be a common suffix
(or prefix) of $\sigma(a)$, $\sigma(b)$ and $\sigma(c)$.

\end{thm}

\bigskip
The following examples illustrate the decomposition after a cyclic
operator:

\begin{example}\label{ex2.1}{}~

\vskip -20 pt{}~
\begin{enumerate}
  \item $\sigma=(ab, acb, ac^2)$. It can be checked that
$\sigma$ is an indecomposable substitution.
 However we have, by a cyclic operation,
$I_{a^{-1}}\circ\sigma=(ba,cba,c^2a)=(b,c,a)\circ
\circ (ac,b,c)\circ (a,b,bc)^2\circ (a,ba,c)$,
which is clearly simple.
\item $\sigma=(acbab, acbaccab, accaccab)$ is indecomposable.
It is easy to check that
$I_{b}\circ\sigma$ is simple.  Notice that in this example,
$I_{a^{-1}}\circ\sigma$ is not simple, though it a substitution.
\end{enumerate}
\end{example}

\begin{rmk}\label{2.99}
  We point out that Theorem \ref{thm1.1} is an easy consequence of
  Theorem \ref{thm2.1}.
\end{rmk}

In fact, let $(w_1,w_2)$ be an invertible substitution,
$w_1,w_2\in\{a,b\}^*$. Then it is obvious that $(w_1,w_2,c)$ is an
invertible substitution over $\{a,b,c\}$. Applying Theorem
\ref{thm2.1}, and noticing that $w$ in the theorem must be empty,
Theorem \ref{thm1.1} follows.

On the other hand, since $I_{b^{-1}}\circ (ba,b)=(ab,b)$, Theorem
\ref{thm1.1} can be restated as the same as Theorem \ref{thm2.1}:

Any invertible substitution over $\{a,b\}$, having been transformed by
some cyclic operator, can be expressed as a finite composition of
the Fibonacci $(ab,b)$ and the permutation $(b,a)$.

\begin{rmk}\label{kj}
It should be noticed that Theorem \ref{thm2.1} can not be extended
directly to the case of an alphabet of 4 letters or more.
\end{rmk}
It suffices to see the following simple example of substitution on
$\{a,b,c,d\}$: $\sigma=(ab,acb,acc,d)$. It is indecomposable, but no inner
automorphism ``$I_z$" can be applied.

\bigskip
Applying the above theorem to the substitution matrices we shall get
the following theorem which explains the Remark \ref{rmk1.2}:
\begin{thm}\label{thm2.2}
  Let $M$ be a $3\times 3$-matrix of non-negative integral
  coefficients.  $M$ is the substitution matrix of some invertible
  substitution if and only if it is a finite product of non-negative
  elementary matrices.
\end{thm}

\section{Proofs}

\indent The theorems will be proved by several lemmas.

We shall use the following notations. The symbol $``+"$(resp.  $``-"$)
will represent various non-empty positive words (resp.  non-empty
negative words), which are the elements of $\langle a, b, c\rangle$
(resp. $\langle a^{-1}, b^{-1}, c^{-1}\rangle$). We also use the
symbols like $``+-+"$ to represent various types of mixed words. As an
example, $u=+-+$ means that $u=u_1u_2^{-1}u_3$, where $u_j\in A^*$ and
$u_1u_2^{-1}u_3$ is the reduced expression
(i.e. $|u_1u_2^{-1}u_3|=|u_1|+|u_2|+|u_3|)$. The meanings of
``$u=+\cdots$" ,``$u=\cdots -$", ``$u=+\cdots -$" \textit{etc}. are
now clear.

Recall that by Nielsen's cancellation method, a morphism
$\sigma=(w_1,w_2,w_3)$ on $\Gamma_A$ is invertible if and only if
$(w_1,w_2,w_3)$ can be carried by finitely many ``elementary Nielsen
transformations" into the identity, and each of these elementary
Nielsen transformations satisfies some ``cancellation" conditions (see
\cite{MKS}). The following lemma is a simple version of Nielsen's
cancellation procedure that we use in the proofs.  For the details, we
refer the reader to \cite{MKS,Ni1,Ni2}.

\begin{lem}\label{lem3.1}
 Let $\sigma=(w_1,w_2,w_3)\in \Aut(\Gamma_A)$.  There exist $k\ge 0$
 and $\tau_1,\cdots,\tau_k\in
\{\pi_1,\pi_2,\phi_l,\phi_r,\iota_1\}$
 such that, denoting
\begin{equation}\label{eqn3.2}
  \sigma_0=\sigma,\ \sigma_i=\sigma_{i-1}\circ\tau_i\ (i=1,\cdots,k),
\end{equation}
\begin{equation}\label{eqn3.3}
 |\sigma_{i}|\le|\sigma_{i-1}|\ (i=1,\cdots, k)
\end{equation}
and that $\sigma_k$ is the identity.
\end{lem}

\bigskip

\noindent \textbf{Notation:} We write the procedure in the above lemma
as $$\sigma\rightarrow\sigma_1\rightarrow\cdots\rightarrow\sigma_k.$$ We
say that each arrow represents a ``Nielsen's cancellation" and that
$\sigma$ is cancelled to $\sigma_k$.
\bigskip

\begin{definition}\label{def3.1}
 We say that a non-trivial substitution $\sigma=(w_1, w_2, w_3)$ is
mixed if it satisfies $w_iw_k^{-1}w_j=+-+$\  $\forall i, j, k\in\{1, 2,
3\}$, $i\not= k$, $j\not= k$.  \end{definition}

The following lemma is proved in \cite{WZ}.
\begin{lem}\label{lem3.2}
Any mixed substitution is non-invertible.
\end{lem}

\begin{cor}\label{cor3.1}
  Let $\sigma=(w_1,w_2,w_3)$ be an invertible substitution, suppose
that for all $i\not=j$, $w_i$ is neither a prefix nor a suffix of
$w_j$.  There exist $\pi\in {\P}$, non-empty words $u$, $v$, $x$, $y$
such that either $\sigma\circ \pi=(ux, uv, yv)$ or
$\sigma\circ\pi=(uxv, uv, y)$.
\end{cor}
\begin{proof}
  By the above lemma, $\sigma $ is not mixed.  Then by the assumption
  of the lemma we have that, there exist $i,j,k \ (i\not=j, j\not=k)$
  such that $w_i^{-1}w_jw_k^{-1}=-$.

Ignoring some permutation, we can only consider the following two
possibilities:

1. $w_1^{-1}w_2w_3^{-1}=-$ ; 2. $w_1^{-1}w_2w_1^{-1}=-$ .

Case 1: We get $w_2=uv$ and $w_1=ux,$ $w_3=yv$ .

Case 2: We must have $w_2=uv$ and $w_1=uf=gv$ . It's clear that
$|u|\neq |g|$ , otherwise $u=g$ then $w_1=w_2$ which contradicts the
invertibility.

\indent Subcase 2.1: $|g|>|u|$ , then by $w_1=uf=gv$ , we have $g=ux$
and $w_1=uxv$.  Then $\sigma =(uxv,uv,y)$ with $y=w_3$ .

\indent Subcase 2.2: $|g|<|u|$ , then still by $w_1=uf=gv$ we have

$u=gx$ and $v=xf$ , thus $w_1=uf$ , $w_2=uxf$ , $\sigma =(uf,uxf,y)$.
\end{proof}

\begin{lem}\label{lem3.3} Let $w,v\in\Gamma_A$ and $\sigma$ be a
morphism on $\Gamma_A$. We have
$$I_w\circ I_v=I_{wv},\ (I_w)^{-1}=I_{w^{-1}},\ \sigma\circ
I_w=I_{\sigma(w)}\circ\sigma.$$
\end{lem}
\begin{proof} By direct verifications. Let us show the last equality:

Let $u\in \Gamma_A$. We have $\sigma\circ I_w(u)=\sigma(wuw^{-1})
=\sigma(w)\sigma(u)\sigma(w)^{-1}=I_{\sigma(w)}(\sigma(u))$.
\end{proof}

\bigskip
Let $w=\alpha_1\alpha_2\cdots\alpha_{|w|}\in \Gamma_A$ ($\alpha_i\in
\A$). We shall use the following notations:
\begin{equation}\label{hq}
h_i(w)=\alpha_i\ (i=1,2,\cdots,|w|),\ h_{\infty}(w)=h_{|w|}(w).
\end{equation}

\begin{lem}\label{lem3.8}
Let $u\in\Gamma_A$ and $x, y\in A^*$ be non-empty.  Suppose that
\begin{equation}\label{uxyu}
  ux=yu
\end{equation}
Then we have either $u\in A^*$ or $u^{-1}\in A^*$.
\end{lem}
\begin{proof}

  First assume that $h_1(u)=+$. We will show that $u\in A^*$.

  Otherwise, we may write $u=u_1\beta^{-1}u_2$ (the reduced
  expression) where $u_1\in A^*,\ u_1\not=\epsilon$, $\beta\in
  A=\{a,b,c\}$ and $u_2\in\Gamma_A$.

  Then by (\ref{uxyu}), we have
\begin{equation}\label{u1}
  u_1\beta^{-1}u_2x=yu_1\beta^{-1}u_2
\end{equation}
Notice that the right hand of the above equation is a reduced
expression, and $$h_{|u_1|+1}(yu_1\beta^{-1}u_2)=+,$$ which,
together with (\ref{u1}), implies
\begin{equation}\label{ttt}
  \beta^{-1}u_2x=+\cdots.
\end{equation}

Since $\beta^{-1}u_2$ is a reduced expression and $x\in A^*$,
(\ref{ttt}) implies that $\beta^{-1}u_2x=+$.  Returning to (\ref{u1}),
we have that the left hand of (\ref{u1}) is a positive word, which
contradicts the right hand.

For the case $h_1(u)=-\ $, we can show similarly that $u^{-1}\in A^*$.
\end{proof}

\begin{cor}\label{cor3.8}
  Let $w\in\Gamma_A$, $\sigma=(w_1,w_2,w_3)$ and $\tau=(u_1,u_2,u_3)$
  be two substitutions.  If $I_w\circ\sigma=\tau$, we have
  $|w_i|=|u_i|$ ($i=1,2,3$).  In particular $|\sigma|=|\tau|$.
\end{cor}
\begin{proof}
  We have $w^{-1}w_iw=u_i\ (i=1,2,3)$. By the above lemma, we have
  $w\in A^*$ or $w^{-1}\in A^*$.  If $w\in A^*$, we have $w_iw=wu_i$
  which implies $|w_i|=|u_i|$; if $w'=w^{-1}\in A^*$, we have
  $w'w_i=u_iw'$.
\end{proof}

\begin{lem}\label{t1}
Let $u,v,x,y\in A^*$. We have

(1) \ If $u$ is not a suffix of $ux$, then $u$ is not a suffix of
$ux^m$ ($m\ge 1$).

(2) \ If $v$ is not a prefix of $yv$, then $v$ is not a prefix of
$y^mv$ ($m\ge 1$).
\end{lem}
\begin{proof}
  We only show that $u\rhd ux^m$ implies $u\rhd ux$ by induction on
  the length of $u$:

  If $|u|\le |x|$, the conclusion is obvious. Suppose now
  $|u|>|x|$. Since $u\rhd ux^m$, we can write $u=u'x$ for some $u'\in
  A^*$, then $u'x\rhd u'xx^m$, that is, $u'\rhd u'x^m$.  Since
  $|u'|<|u|$, by the hypothesis of induction we obtain $u'\rhd u'x$,
  which implies $u\rhd ux$.
\end{proof}

The following lemmas study the ``cancellation properties" between
special words, where we shall use $h_1(w)$ and $h_{\infty}(w)$ to
indicate the first and the last character of $w$ as defined in
(\ref{hq}).

\begin{lem}\label{lem3.4}
Let $u$, $v$, $x$ and $y$ be non-empty words satisfying
\begin{eqnarray}
& & \ h_1(u)\not=h_1(y)\textrm{ and } h_{\infty}(x)\not=h_{\infty}(v);
\label{e1}\\ & & \ x \ (\textrm{resp. } v)\textrm{ is not a prefix of
} v\ (\textrm{resp. } x); \label{e2} \\ & & \ y \ (\textrm{resp. }
u)\textrm{ is not a suffix of }u\ (\textrm{resp. } y)\label{e3}.
\end{eqnarray}

Let $w_1=ux,w_2=uv,w_3=yv$ and consider the following mixed word:
\begin{equation}\label{eqn3.7}
  w=w_{i_0}^\epsilon w_{i_1}^{-\epsilon}\cdots
  w_{i_k}^{(-1)^k\epsilon}
\end{equation}
where $k\ge 0$, $\epsilon\in\{+1,-1\}$, $i_m\in\{1,2,3\}\
(m=0,1,\cdots,k)$ and $i_m\not=i_{m+1}\ (m=0,1,\cdots,k-1)$.  Then we
have
\begin{equation}\label{eqn3.8}
  h_1(w)=h_1(w_{i_0}^\epsilon),\ h_{\infty}(w)=
  h_{\infty}(w_{i_k}^{(-1)^k\epsilon})
\end{equation}
and
\begin{equation}\label{eqn3.88}
|w|\ge 2.
\end{equation}
(Notice that (\ref{eqn3.7}) is in general not a reduced
expression, this means that there are possible cancellations
between words).
\end{lem}
\begin{proof}
 First notice the following obvious facts:

 Let $u_1, u_2\in {\Gamma_{A}}^*$ be two (mixed) words in reduced
 expression.  $|u_1u_2|<|u_1|+|u_2|$ (i.e. some cancellation
 occurs) if and only if $h_{\infty}(u_1)=(h_1(u_2))^{-1}$.

If $h_{\infty}(u_1)\not=(h_1(u_2))^{-1}$, then $v=u_1u_2$ is a reduced
expression, thus $h_1(v)=h_1(u_1)$ and $h_{\infty}(v) =h_{\infty}
(u_2)$.

  We only show $h_1(w)=h_1(w_{i_0}^\epsilon)$. The proof for
$h_{\infty}(w)=h_{\infty}(w_{i_k}^{(-1)^k\epsilon})$ is similar.

We consider only the case $\epsilon=+1$, the case $\epsilon=-1$ is
similar.

When $k\le 3$, we can check the conclusion directly by enumerating all
possible cases. Now we will show it for $k\ge 4$ by induction.

We have to consider carefully several cases. Keep in mind that in the
following argument, we shall use frequently the condition
(\ref{e1})--(\ref{e3}) and the facts mentioned at the beginning of the
proof.

\bigskip
\noindent\textsl{Case 1.}  $w_{i_1}=uv$

\indent \textsl{Subcase 1.1.} $w_{i_0}=ux$:

By the hypothesis of induction $h_1$
$(w_{i_1}^{-1}\cdots)$
$=h_1(v^{-1})$,
hence $h_1(w)=h_1(uxv^{-1}\cdots)$
$=h_1(u)=h_1(w_{i_0})$.

\indent \textsl{Subcase 1.2.} $w_{i_0}=yv$

\indent \textsl{Sub-Subcase 1.2.1.} $w_{i_2}=yv$:

Since $h_1(w_{i_2}\cdots)=h_1(y)$, we have
$h_1(w)=h_1(yv(uv)^{-1}y\cdots)=h_1(y)=h_1(w_{i_0}).$

\indent \textsl{Sub-Subcase 1.2.2.} $w_{i_2}=ux$:

 In this case, either $w_{i_3}=uv$ or $w_{i_3}=yv$, we always have
 $h_1(w_{i_3}^{-1}\cdots)=h_1(v^{-1})$. Thus $h_1(w)=
 h_1(yv(uv)^{-1}uxv^{-1} \cdots) =h_1(yxv^{-1}\cdots)
 =h_1(y)=h_1(w_{i_0})$.

\bigskip
\noindent\textsl{Case 2.} $w_{i_1}=yv$

\indent \textsl{Subcase 2.1.} $w_{i_0}=ux$:

$h_1(w_{i_1}^{-1}\cdots)=h_1(v^{-1})$ and
$h_1(w)=h_1(uxv^{-1}\cdots)=h_1(u)=h_1(w_{i_0})$.

\indent \textsl{Subcase 2.2.} $w_{i_0}=uv$:

In this case, either $w_{i_2}=ux$ or $w_{i_2}=uv$, we always have
 $h_1(w_{i_2}\cdots)=h_1(u)$. Thus $h_1(w)=h_1(uv(yv)^{-1}u\cdots)=
 h_1(uy^{-1}u\cdots) =h_1(u)=h_1(w_{i_0})$.

\bigskip
\noindent\textsl{Case 3.}  $w_{i_1}=ux$

\indent \textsl{Subcase 3.1.} $w_{i_0}=uv$:

$h_1(w_1^{-1}\cdots)=h_1(x^{-1})$ and
$h_1(w)=h_1(uvx^{-1}\cdots)=h_1(u)=h_1(w_{i_0})$.

\indent \textsl{Subcase 3.2.} $w_{i_0}=yv$:

$h_1(w_{i_1}^{-1}\cdots)=h_1(x^{-1})$ and
$h_1(w)=h_1(yvx^{-1}\cdots)=h_1(y)=h_1(w_{i_0})$.
\end{proof}

The above lemma will be used to study substitutions of the form
$(ux,uv,yv)$.  The following lemma is the version for $(uxv,uv,y)$.

\begin{lem}\label{t2}
Let $u$, $v$, $x$ and $y$ be non-empty words satisfying
\begin{eqnarray}
& & \ h_1(u)\not=h_1(y)\textrm{ and } h_{\infty}(v)\not=h_{\infty}(y);
\nonumber\\ & & \ v\textrm{ is not a prefix of }xv;\nonumber \\ & & \
u\textrm{ is not a suffix of }ux. \nonumber
\end{eqnarray}

Let $w_1=uxv,w_2=uv,w_3=y$ and consider the mixed words of the form
(\ref{eqn3.7}). Then

(i)\ (\ref{eqn3.8}) holds;

(ii)\ If $|y|>1$, (\ref{eqn3.88}) always holds;

(iii)\ If $|y|=1$, (\ref{eqn3.88}) holds except when $k=0$, and
$w_{i_0}=w_3$.

\end{lem}
\begin{proof}

(ii) and (iii) are consequences of (i). To show (i), as in the above
lemma, we only show $h_1(w)=h_1(w_{i_0}^\epsilon)$ where
$\epsilon=+1$.

  The cases $k=0,1$ are trivial. Let us show the lemma for $k\ge 2$ by
  induction.

\bigskip
\noindent\textsl{Case 1.} $w_{i_0}=uxv$

\textsl{Subcase 1.1.} For all $j$, $w_{i_j}\not=y$:

Then $w_{i_1}=uv, w_{i_2}=uxv, \cdots$.

\textsl{Subsubcase 1.1.1.} $k=2m-1$:

In this case, $w=\left((uxv)(uv)^{-1}\right)^m=ux^mu^{-1}$.  Then by
Lemma \ref{t1} we know that $w=+-$ which implies
$h_1(w)=h_1(u)=h_1(w_{i_0})$.

\textsl{Subsubcase 1.1.2.} $k=2m$:

$w=ux^{m+1}v$, the conclusion follows.

\textsl{Subcase 1.2.}  Let $w_{i_j}=y$, where $j=\min\{m; w_{i_m}=y\}$

\textsl{Subsubcase 1.2.1.} $j=2m-1$:

We have $w=\left((uxv)(uv)^{-1}\right)^my\cdots=ux^mu^{-1}y\cdots$,
thus $h_1(w)=h_1(u)=h_1(w_{i_0})$.

\textsl{Subsubcase 1.2.2.} $j=2m$:

We have $w=\left((uxv)(uv)^{-1}\right)^m(uxv)y^{-1}\cdots
=ux^{m+1}vy^{-1}\cdots$ then the conclusion follows.

\bigskip
\noindent\textsl{Case 2.}  $w_{i_0}=uv$. We use a similar proof as the
case 1.

\bigskip
\noindent\textsl{Case 3.}  $w_{i_0}=y$:

For any subcases, it is checked that $h_1(w_{i_1}^{-1}\cdots
w_{i_k}^{(-1)^k)})=h_1(v^{-1})$, thus

\hfill $h_1(w)=h_1(yv^{-1}\cdots)=h_1(y)=h_1(w_{i_0}).$
\end{proof}

\bigskip
\begin{lem}\label{lem3.5}
Under the notations and conditions of Lemma \ref{lem3.4} (resp. Lemma
 \ref{t2}), let $\sigma=(w_1,w_2,w_3)$, then $\sigma$ is not
 invertible.
\end{lem}

\begin{proof}   Were $\sigma$ invertible, we would have
a Nielsen's cancellation procedure as follows
\begin{equation}\label{tmp}
  \sigma\rightarrow\cdots\rightarrow
(w_1^{(k)},w_2^{(k)},w_3^{(k)})\rightarrow(w_1^{(k+1)},
w_2^{(k+1)},w_3^{(k+1)}) \cdots\rightarrow I
\end{equation}
where
\begin{equation}\label{tmp11}
  w_j^{(k+1)}\not=\varepsilon \ (j=1,2,3).
\end{equation}

We will show that:

\bigskip
\noindent\textsl{Claim:} $w_i^{(k)}$ ($i=1,2,3,\ k=1,2,\cdots$) can be
written in the form (\ref{eqn3.7}).

\bigskip
Then the previous lemma implies that $|w_i^{(k)}|\ge 2$ which is in
contradiction with (\ref{tmp}).

Now let us prove the claim by induction.

Suppose that $w_1^{(k)},w_2^{(k)},w_3^{(k)}$ are of the form
(\ref{eqn3.7}).  Then $$(w_1^{(k+1)},w_2^{(k+1)},w_3^{(k+1)})=
(w_1^{(k)},w_2^{(k)},w_3^{(k)})\circ \tau$$ where
$\tau\in\{\pi_1,\pi_2,\phi_l,\phi_r\,\iota_1\}$.  We shall show that
$w_1^{(k+1)},w_2^{(k+1)},w_3^{(k+1)}$ are still of the form
(\ref{eqn3.7}).

If $\tau\in\{\pi_1,\pi_2,\iota_1\}$, the conclusion is trivial.
So we consider only the case $\tau=\phi_r$. The case for
 $\tau=\phi_l$ can be treated similarly.

Now $w_1^{(k+1)}=w_1^{(k)}w_2^{(k)}$ (and $w_2^{(k+1)}=
w_2^{(k)},w_3^{(k+1)}=w_3^{(k)}$), then by the cancellation condition,
we must have $|w_1^{(k)}w_2^{(k)}|\le |w_1^{(k)}|$ which implies
\begin{equation}\label{temp1}
 h_{\infty}(w_1^{(k)})h_1(w_2^{(k)})=\varepsilon.
\end{equation}
Without loss of generality we suppose that $h_{\infty}(w_1^{(k)})=-$,
that is (by the above lemma and the assumption of induction),
$w_1^{(k)}=\cdots w_{\lambda_1}w_{\lambda_2}^{-1}$.  Then by the
cancellation condition (\ref{temp1}) we must have $w_2^{(k)}=+$, that
is, $w_2^{(k)}=w_{\mu_1}w_{\mu_2}^{-1}\cdots$.  Consequently
$w_1^{(k+1)}=\cdots
w_{\lambda_1}w_{\lambda_2}^{-1}w_{\mu_1}w_{\mu_2}^{-1} \cdots$.  If
$\lambda_2\not=\mu_1$, the conclusion for $w_1^{(k+1)}$ is proved; if
$\lambda_2=\mu_1$, by cancelling $w_{\lambda_2}^{-1}w_{\mu_1}$ in the
above expression and by induction, using the fact that
$w_1^{(k+1)}\not=\varepsilon$ (by \ref{tmp11}),
we can show easily that $w_1^{(k+1)}$
is still of the form (\ref{eqn3.7}).  The claim is thus proved.
\end{proof}

The following remarks explain the original idea of the above proof:
\begin{rmk}\label{addrmk2}
  There is an essential difference between the proof of the above
  lemmas and that of Theorem 2.1 in \emph{\cite{WZ}} : here
  $w_1,w_2,w_3$ are not mixed.  The essential point is that, by the
  assumption of the lemma, we can control the possible cancellations
  such as $w_1^{-1}w_2w_3^{-1}\not=-+-$.

\begin{rmk}\label{addrmk1}
  An alternative proof can be done as follows, we omit the details
  which are not trivial: starting from $\sigma$, consider all possible
  Nielsen's cancellations procedures, we shall get a finite tree of
  the procedures (which contains several periodic circle), no branch
  of the tree is reduced to the identity.
\end{rmk}
\bigskip

\end{rmk}
\begin{prp}\label{prp3.1}
Suppose that $\sigma=(w_1,w_2,w_3)$ is a non-simple indecomposable
substitution.  Then we have either $h_1(w_1)=h_1(w_2)=h_1(w_3)$ or
$h_{\infty}(w_1) =h_{\infty}(w_2)=h_{\infty}(w_3)$.  In other words,
$w_1,w_2,w_3$ must have a common non-empty prefix or suffix.
\end{prp}

\begin{proof}
Since $\sigma$ is indecomposable and invertible, $\sigma$ satisfies
the condition in Corollary \ref{cor3.1}.  Hence we can suppose that

$
\begin{array}{lr}
  \textrm{Case 1. } \sigma=(w_1,w_2,w_3)=(ux, uv, yv) & \\
  \textrm{Case 2. } \sigma=(uxv, uv, y) &
\end{array}
$

Consider the case 1.

If the conclusion of the proposition were not true, we would have that
$(ux, uv, yv)$ satisfies the condition of Lemma \ref{lem3.4} and
\ref{lem3.5}, which imply that $\sigma$ is not invertible, a
contradiction.

For the case 2, the proof is the same as the Case 1 except that we
apply Lemma \ref{t2} in stead of Lemma \ref{lem3.4}.
\end{proof}

\begin{lem}\label{prelem}
 Let $z,u,v,x,y\in A^*$. Assume that $u,v,x,y$ satisfy (\ref{e1}),
 (\ref{e2}) and (\ref{e3}).  Denote $w_1=ux,w_2=uv,w_3=yv$,
$$
\begin{array}{rl}

\mathcal{N}_1=&\{ zw_{i_0}w_{i_1}^{-1}\cdots w_{i_{2k-1}}^{-1}z^{-1};
\ k\ge 1, i_m=1,2,3,\\
&(m=0,\cdots,2k-1),\ i_m\not=i_{m+1}\ (m=0,\cdots,2k-2)\} \\

\mathcal{N}_2=& \{ zw_{i_0}w_{i_1}^{-1}\cdots w_{i_{2k}};
\ k\ge 1, i_m=1,2,3, \\
&(m=0,\cdots,2k),\ i_m\not=i_{m+1}\
(m=0,\cdots,2k-1)\} \\

\mathcal{N}_3=&  \{ w_{i_0}^{-1}w_{i_1}\cdots w_{i_{2k-1}};
\ k\ge 1, i_m=1,2,3, \\
&(m=0,\cdots,2k-1),\ i_m\not=i_{m+1}\
(m=0,\cdots,2k-2)\} \\

\mathcal{N}_4=&\{ w_{i_0}^{-1}w_{i_1}\cdots
w_{i_{2k-1}}w_{i_{2k}}^{-1}z^{-1};
\ k\ge 1, i_m=1,2,3, \\
&(m=0,\cdots,2k),\ i_m\not=i_{m+1}\ (m=0,\cdots,2k-1)\}\\

\mathcal{N}=& \mathcal{N}_1\bigcup\mathcal{N}_2\bigcup\mathcal{N}_3
\bigcup\mathcal{N}_4
\end{array}$$

We have
$$
\begin{array}{llr}
(1) &\textrm{(\ref{eqn3.8}) holds;} \\ (2) &\textrm{If
}(zw_1,zw_2,zw_3) \textrm{ is cancelled by a Nielsen's procedure to }
(U_1,U_2,U_3), \textrm{ then }\\ & U_1,U_2,U_3 \in \mathcal{N};\\ (3)
&(zux,zuv,zyv)\textrm{ is not invertible.}\\ (4)
&(uxz,uvz,yvz)\textrm{ is not invertible.}
\end{array}$$

\end{lem}

\begin{proof} (1) is proved in Lemma \ref{lem3.4}.
(3) is a consequence of (2). (4) can be proved similarly (by using
similar definitions of $\mathcal{N}_k$, or simply by the fact that
(4) differs from (3) by an inner automorphism).  Let us prove (2) by
induction.

Suppose that $\sigma\rightarrow\cdots\rightarrow(V_1,V_2,V_3)
\rightarrow(U_1,U_2,U_3)$ where $V_1,V_2,V_3\in \mathcal{N}$, we show
that $U_1,U_2,U_3\in \mathcal{N}$. With a simple observation of
appearances of $z$ and $z^{-1}$ in the expressions above, by the
cancellation condition, the conclusion can be proved by using
analogous arguments as in the proof of Lemma \ref{lem3.4}.
\end{proof}

\begin{rmk}\label{t3}
  When $x,y,u,v,w_1,w_2,w_3$ are given as in Lemma \ref{t2}, similar
  conclusions of the above lemma hold. The proof is the same.
\end{rmk}


\begin{lem}\label{lem3.6}
Suppose that $\sigma=(zu_1,zu_2,zu_3)$ (\emph{resp.}
$\sigma=(u_1z,u_2z,u_3z)$) is a non-simple indecomposable substitution where
$z$ is a non-empty word. Suppose that $u_1,u_2,u_3$ have no common
prefix and no common suffix.  Then there exist a non-trivial
invertible substitution $\sigma'$, a permutation $\pi\in\P$ and
a Fibonacci $f\in\F$, such that
\begin{equation}\label{eqn3.9}
  I_{z^{-1}}\circ\sigma=\sigma'\circ f\circ \pi \qquad
  (\emph{resp. }I_{z}\circ\sigma=\sigma'\circ f\circ \pi).
\end{equation}
\end{lem}
\begin{proof}
Let $w_i=zu_i$ ($i=1,2,3$). First, from the fact that $\sigma$
is non-simple indecomposable we have
\begin{equation}\label{z1}
  w_i\textrm{ is neither a prefix nor a suffix of } w_j\ (\forall
  i,j,\ i\not=j).
\end{equation}

By Lemma \ref{lem3.2}, $(w_1,w_2,w_3)$ is not mixed. It turns out that
$(u_1,u_2,u_3)$ is not mixed. By Corollary \ref{cor3.1}, ignoring some
permutation, we may write
\begin{equation}\label{only}
  (u_1,u_2,u_3)=(ux,uv,yv)
\end{equation}
\begin{equation}\label{only1}
  \textrm{ or \ }(u_1,u_2,u_3)=(uxv,uv,y).
\end{equation}

We only prove the lemma for the case (\ref{only}). The proof for the
case (\ref{only1}) is the same, because of Lemma \ref{t2} and Remark
\ref{t3}.

Now we have $\sigma=(zux,zuv,zyv)$. Then
$I_{z^{-1}}\circ\sigma=(uxz,uvz,yvz)$.  We claim that there exist a non-trivial
substitution $\sigma'$ and a Fibonacci substitution $f$ such
that
\begin{equation}\label{tr}
  (uxz,uvz,yvz)=\sigma'\circ f
\end{equation}
In fact, if (\ref{tr}) fails to hold, denoting
$(uxz,uvz,yvz)=(w_1',w_2',w_3')$, we shall have
\begin{equation}\label{z2}
w_i'\textrm{ is neither a prefix nor a suffix of } w_j'\ (\forall
  i,j,\ i\not=j).
\end{equation}
then it follows from (\ref{z1}) and (\ref{z2}) that $u,v,x,y$ satisfy
the condition of Lemma \ref{prelem}, which implies that $\sigma$ is
not invertible, a contradiction.  Hence $I_{z^{-1}}\circ \sigma$ can be
decomposed as in (\ref{tr}) and the lemma follows.
\end{proof}


 \begin{lem}\label{lem3.7} Let $\sigma$ be a non-simple indecomposable
substitution.  Then there exists $z\in\Gamma_A$ such that
$I_z\circ\sigma$ is decomposable. That is, there exist
non-trivial invertible
substitutions $\sigma_1$ and $\sigma_2$ such that
\begin{equation}\label{eqn3.10}
  I_z\circ\sigma=\sigma_1\circ \sigma_2
\end{equation}
Furthermore,
$|\sigma_i|<|\sigma| \ (i=1,2)$.
\end{lem}

\begin{proof}
Let $\sigma=(w_1,w_2,w_3)$. Since it is indecomposable, for any $i\not=j$,
$w_i$ is not a prefix (resp. suffix) of $w_j$. Hence we may write
$$(w_1,w_2,w_3)=(u_1w,u_2w,u_3w)$$ where $w$ is the maximal common
suffix of $w_1,w_2,w_3$ (if any) and $u_1,u_2,u_3$ are non-empty
words.

Then we have:
$$I_{w}\circ\sigma=(wu_1,wu_2,wu_3)$$ where $u_1,u_2,u_3$ have no
common suffix.

If $(wu_1,wu_2,wu_3)$ is not indecomposable, the lemma is proved.

Now suppose that $(wu_1,wu_2,wu_3)$ is indecomposable.  Since $wu_1,wu_2,wu_3$
have no common suffix, by Proposition \ref{prp3.1}, we can write
$$
  I_{w}\circ\sigma=(wu_1,wu_2,wu_3)=(zu_1',zu_2',zu_3')
$$
where $z\not=\varepsilon$ is the maximal common prefix of $wu_1,wu_2$,
and $wu_3$. In this case, $u_1',u_2',u_3'$ have no common prefix and
no common suffix. Then by Lemma \ref{lem3.6}, $I_{z^{-1}}\circ
I_{w}\circ\sigma$ is decomposable. That is (by Lemma
\ref{lem3.3}) $I_{z^{-1}w}\circ\sigma$ is decomposable.

The last statement of the lemma is an easy consequence of Corollary
\ref{cor3.8}.
\end{proof}

Now we are ready to prove our theorems.

\bigskip
\noindent \textbf{Proof of Theorem \ref{thm2.1}:}

\bigskip

 First we prove Theorem \ref{thm2.1} (1). Let $\sigma$ be an
 invertible substitution. Then there exist $k\ge 1$ and indecomposable
 substitutions $\sigma_1, \cdots, \sigma_k$ such that
 $\sigma=\sigma_1\circ\cdots\sigma_k$.

If some $\sigma_i$ is non-simple, then by Lemma \ref{lem3.7}, there
exist $z\in\Gamma_A$, invertible substitutions $\sigma_i'$ and
$\sigma_i''$ such that
$\sigma_i=I_{z^{-1}}\circ\sigma_i'\circ\sigma_i''$ and that
\begin{equation}\label{zzz5}
  |\sigma_i'|<|\sigma_i|,\ |\sigma_i''|<|\sigma_i|.
\end{equation}
Then we repeat such decomposition for $\sigma_i'$ (resp.
$\sigma_i''$) and so on. By (\ref{zzz5}), such decomposition will
terminate after finite times, that is, finally every factor will be a
simple substitution. Hence we can write
$$\sigma=I_{w_1}\circ\tau_1\circ I_{w_2}\circ\tau_2\cdots\circ
I_{w_n}\circ\tau_n$$ where $\tau_i\ (i=1,\cdots,n)$ is a simple
substitution and $w_i\in\Gamma_A$ (put $w_i=\varepsilon$ if
necessary).

By Lemma \ref{lem3.3} (and a simple induction), we may write
$$\sigma=I_{w}\circ\tau_1\circ\tau_2\cdots\circ\tau_n$$ where $w\in
\Gamma_A$.

Let $\tau=\tau_1\circ\tau_2\cdots\circ\tau_n$. It is clear that $\tau$
is a simple substitution.

Finally, since $I_{w^{-1}}\circ\sigma=\tau$, Lemma \ref{lem3.8}
implies that $w\in A^*$ or $w^{-1}\in A^*$. Theorem \ref{thm2.1}(1) is
thus proved.

\bigskip
For Theorem \ref{thm2.1}(2), denoting
$$|\sigma|_{min}=\min\{|\sigma(a)|, |\sigma(b)|,|\sigma(c)|\},$$ it is
easy to see that it is equivalent to prove the following claim:

\bigskip
\noindent \textsl{Claim:} We can choose $w$ in (1) such that $|w|\le
|\sigma|_{min}$.

\bigskip
Let us prove the claim by induction on $k=|\sigma|$.

If $k\le 4$ the claim can be verified simply by enumerating all cases.

Suppose that the claim is true for $k\le n$ (more clearly, for any
substitution $\sigma'$ of length less that $n$) and that
$|\sigma|=n+1$. By the conclusion (1) of the theorem, there exists
$w\in A^*$ or $w^{-1}\in A^*$ such that $I_{w}\circ \sigma$ is
simple. To be specific, we can suppose that $w^{-1}\in A^*$,
$\sigma=(w_1,w_2,w_3)$ and that $|w_1|=|\sigma|_{min}$.

If $|w|\le |\sigma|_{min}$, nothing need to prove. Suppose
$|w|>|\sigma|_{min}=|w_1|$. It is easy to see that $w_1$ is a common
prefix of $w_1, w_2, w_3$. Let then $w_2=w_1w_2'$. We have
$(w_1,w_2,w_3)=(w_1,w_2',w_2) \circ (a,ab,c)$, that is,
\begin{equation}\label{today1}
  \sigma=\sigma'\circ f.
\end{equation}
where $\sigma'=(w_1,w_2',w_2)$ is obviously an invertible substitution
and $f=(a,ab,c)$ is a Fibonacci substitution.

It is trivial that $|\sigma'|<|\sigma|=n+1$ and that
$|\sigma'|_{min}\le |\sigma|_{min}$ , hence by the hypothesis of
induction we have the following fact:

There exists $z\in A^*$ (or $z^{-1}\in A^*$) such that $|z|\le
|\sigma'|_{min}$ and that
\begin{equation}\label{today2}
  g:=I_z\circ \sigma' \textrm{ is simple.}
\end{equation}

Then it follows from (\ref{today1}) and (\ref{today2}) that
$$I_z\circ \sigma=I_z\circ \sigma'\circ f=g\circ f.$$

Since $|z|\le|\sigma|_{min} $ and $g\circ f$ is (by definition)
simple, the conclusion is proved.  \hfill$\Box$

\bigskip
\noindent \textbf{Proof of Theorem \ref{thm2.2}:}

\bigskip
We check easily that if $\sigma=I_w\circ \tau$, then
$M_{\sigma}=M_{\tau}$. We have also $M_{\sigma\circ\tau}=M_{\sigma}
M_{\tau}$.  The theorem is then a direct consequence of Theorem
\ref{thm2.1}.\hfill$\Box$

\bigskip
\noindent \textbf{Acknowledgement:} \ The authors would like to thank
Professors S. Ito, J. Peyri\`ere, and Z.-Y. Wen for helpful
discussions.


\end{document}